\documentclass[12pt]{amsart}

\usepackage{amsthm}
\usepackage{amsmath, mathrsfs}
\usepackage{amssymb}
\usepackage{hyperref}
\usepackage{enumerate}
\usepackage{latexsym,array}
\usepackage{amsfonts}
\usepackage{shadow}
\usepackage{tikz}
\usepackage{rotating}
\usepackage{tabularx}
\usepackage{graphicx}
\usepackage[section]{placeins}
\usepackage{mathtools}
\usepackage{url}

\newcommand{\bigE}{\mathcal{E}}

\newcommand{\bigK}{\mathcal{K}}

\newcommand{\bigG}{\mathcal{G}}

\newcommand{\bigO}{\mathcal{O}}
\newcommand{\bigT}{\mathcal{T}}
\newcommand{\bigV}{\mathcal{V}}

\newtheorem{Pa}{Paper}[section]
\newtheorem{Tm}[Pa]{{\bf Theorem}}
\newtheorem{La}[Pa]{{\bf Lemma}}
\newtheorem{Dn}[Pa]{{\bf Definition}}
\newtheorem{Cy}[Pa]{{\bf Corollary}}

\newtheorem{Rk}[Pa]{{\bf Remark}}
\newtheorem{Pn}[Pa]{{\bf Proposition}}

\newtheorem{Ex}[Pa]{{\bf Example}}

\tikzset{node distance=2cm, auto}

\date{}

\title[Hyperrigid subsets of Cuntz--Krieger algebras]
{Hyperrigid subsets of Cuntz--Krieger algebras and the property of\\ rigidity at zero}

\author[G. Salomon]{Guy Salomon}
\address{Pure Mathematics Department\\University of Waterloo\\Waterloo\\Ontario N2L 3G1\\Canada}
\email{gsalomon@uwaterloo.ca}

\subjclass[2010]{47L80, 47L55, 47L40, 46L05}
\keywords{Cuntz--Krieger algebra; directed graph; hyperrigidity; C*-envelope}
\thanks{The author was partially supported by the Clore Foundation.}

\begin{document}

\maketitle
\begin{abstract}
A subset $\mathcal{G}$ generating a $C^*$-algebra $A$ is said to be {\em hyperrigid} if for every faithful nondegenerate $*$-representation $A\subseteq B(H)$ and a sequence $\phi_n:B(H) \to B(H)$ of unital completely positive maps, we have that
\[
\lim_{n\to\infty}\phi_n(g)= g~~\text{for all } g\in \mathcal{G} ~~ \implies ~~ 
\lim_{n\to\infty}\phi_n(a)= a~~\text{for all } a\in A
\]
where all convergences are in norm. In this paper, we show that for the Cuntz--Krieger algebra $\mathcal{O}(G)$ associated to a row-finite directed graph $G$ with no isolated vertices, the set of partial isometries $\mathcal{E}=\{S_e:e\in E\}$ is hyperrigid.

In addition, we define and examine a closely related notion: the property of rigidity at $0$. A generating subset $\mathcal{G}$ of a $C^*$-algebra $A$ is said to be {\em rigid at $0$} if for every sequence of contractive positive maps $\varphi_n:A\to \mathbb C$ satisfying $\lim_{n\to \infty}\varphi_n(g)=0$ for every $g\in \mathcal{G}$, we have that $\lim_{n\to \infty}\varphi_n(a)=0$ for every $a\in A$.

We show that, when combined, hyperrigidity and rigidity at $0$ are equivalent to a somewhat stronger notion of hyperrigidity, and we connect this to the unique extension property. This, however, is not the case for the generating set $\mathcal{E}$. More precisely, we show that for any graph $G$, subsets of the Cuntz--Krieger family generating $\mathcal{O}(G)$ are rigid at $0$ if and only if they contain every vertex projection.
\end{abstract}

\section{Introduction}

A {\em directed graph} $G=(V,E,s,r)$ consists of a set $V$ of vertices, a set $E$ of edges, and two maps $s,r : E\to V$, called the {\em source} and the {\em range}; if $v=s(e)$ and $w=r(e)$ we say that $v$ emits $e$ and $w$ receives it. In this paper, we consider only {\em countable} directed graphs, meaning that both the sets $V$ and $E$ are countable.
A directed graph is said to be {\em row-finite} if every vertex receives at most finitely many edges, i.e., $r^{-1}(v)$ is a finite subset of $E$, for all $v\in V$. 

A {\em Cuntz--Krieger $G$-family} $(\mathcal V,\mathcal E)$ of a directed graph $G$ consists of a set of mutually
orthogonal projections $\mathcal V:= \{P_v : v \in V\}$ and a set of partial isometries $\mathcal E:= \{S_e : e \in E\}$ which satisfy the relations:
\begin{enumerate}
\item[(I)] $S_e^*S_e=P_{s(e)}$, for every $e\in E$; 
\item[(TCK)] $\sum_{e\in F} S_e S_e^* \leq P_{v}$ for every finite subset $F \subset r^{-1}(v)$; and
\item[(CK)] $\sum_{r(e)=v} S_e S_e^* = P_v$, for every $v \in V$, with $0<|r^{-1}(v)|<\infty$.
\end{enumerate}
There exists 
a universal $C^*$-algebra $\bigO(G)$ generated by a Cuntz--Krieger $G$-family, which is called the {\em Cuntz--Krieger algebra} of the graph $G$. The original definition of this $C^*$-algebra is due to Cuntz and Krieger \cite{Cuntz-Krieger}; for a comprehensive background on Cuntz--Krieger algebras associated to directed graphs we refer the reader to Raeburn's book \cite{Raeburn}.  

When $G$ is row-finite, it is known \cite{Duncan,Kakariadis,DorOn-Salomon} that inside the Cuntz--Krieger algebra $\bigO(G)$, the Cuntz--Krieger family is not just a generating set, but in fact, a {\em hyperrigid} generating set. 

\begin{Dn}\label{Dn:HR}
Let $A$ be a $C^*$-algebra and $\bigG$ be a generating subset of $A$.
We say that $\bigG$ is {\em hyperrigid} in $A$ if for every faithful nondegenerate $*$-representation $A\subseteq B(H)$ and a sequence of unital completely positive maps $\phi_n:B(H) \to B(H)$ we have that
\[
\lim_{n\to\infty}\phi_n(g)= g~~\text{for all } g\in \bigG ~~ \implies ~~ 
\lim_{n\to\infty}\phi_n(a)= a~~\text{for all } a\in A
\]
where all convergences are in norm.
\end{Dn}

The notion of hyperrigidity was extensively studied during the last decade in various contexts; see for example  \cite{Arv_HR,Davidson-Kennedy,Kats-Hao-Ng,Kennedy-Shalit-HR}  and \cite[Section 7]{Kak-Shalit}. There is one question that naturally arises whenever a set $\bigG$ is hyperrigid in a $C^*$-algebra $A$: what is the smallest subset of $\bigG$ that is already hyperrigid in $A$?

In \cite{Arv_HR}, Arveson showed that $\{1,x,x^2\}$ is hyperrigid in $C([0,1])$ while the smaller generating set $\{1,x\}$ is not. In \cite{Kennedy-Shalit-HR}, Kennedy and Shalit considered the Cuntz algebra $\bigO_I$ associated to a homogeneous ideal $I \triangleleft \mathbb C[z_1,\dots, z_d]$. Without getting into technical details, let us just say that the $C^*$-algebra $\bigO_I$ is naturally generated by a set of $d$ generators $\bigG=\{Z_1,\dots,Z_d\}$. While it can be easily shown that 
the larger generating set $\bigG \cup \bigG^*\bigG$ is always hyperrigid in $\bigO_I$ for every homogeneous ideal $I$, the authors showed in \cite[Theorem 4.12]{Kennedy-Shalit-HR} that hyperrigidity of $\bigG$ itself is equivalent to the well known essential normality conjecture \cite[Problem 2]{ArvDirac02} of Arveson.

In \cite{DorOn-Salomon}, Dor-On and the author showed that $G$ being row-finite does not only imply the hyperrigidity of the Cuntz--Krieger $G$-family inside $\bigO(G)$, but is in fact equivalent to the latter (in the paper, the authors considered the norm-closed algebra generated by the Cuntz--Krieger family --- namely, the tensor-algebra $\bigT_+(G)$ --- but the proof works the same if one replaces $\bigT_+(G)$ with the Cuntz--Krieger family itself; see the proofs of \cite[Theorems 3.5 and 3.9]{DorOn-Salomon}).

It is therefore natural to ask the following question regarding a row-finite graph $G$: in case $\bigO(G)$ is already generated by a subset $\bigG$ of the Cuntz--Krieger family, is this subset hyperrigid as well? For simplicity, if the graph is assumed to contain no isolated vertices, then the subset $\bigE:=\{S_e:e\in E\}$ is a (minimal) subset of the Cuntz--Krieger family generating $\bigO(G)$; must $\bigE$ be hyperrigid in this case?

In this paper, we show that answer is {\em yes}: $\bigE$ is hyperrigid in $\bigO(G)$ if and only if $G$ is a row-finite.
In the proof, we use a generalization to a not-necessarily-unital $C^*$-algebra $A$ of Arveson's characterization for the hyperrigidity of a generating subset $\bigG$: $\bigG$ is hyperrigid in $A$ if and only if for every unital $*$-representation $\pi:A \to B(H)$ there exists a unique unital completely positive extension of $\pi|_{\bigG}$ to $A$, namely $\pi$ itself. We also prove that in a proper sense hyperrigidity is preserved under inductive limits, and use this to show the above result.

The fact that for a row-finite graph $G$, whenever $\bigE$ generates $\bigO(G)$ it must be hyperrigid,  also implies that for row-finite graphs the $C^*$-envelope of $A(\bigE):=\overline{{\rm{alg}}}\{S_e:e\in E\}$ --- the operator algebra generated by $\bigE$ --- is $\bigO(G)$. This is known when the algebra $A(\mathcal E)$ is replaced by the tensor algebra $\bigT_+(G)=\overline{{\rm{alg}}}\{P_v,S_e:v\in V, e\in E\}$; see \cite[Theorem 2.5]{KatsoulisKribs_Graphs}.

In addition, we define and examine a closely related notion: the property of rigidity at $0$. A generating subset $\bigG$ of a $C^*$-algebra $A$ is said to be {\em rigid at $0$} if for every sequence of contractive positive maps $\varphi_n:A\to \mathbb C$ satisfying $\lim_{n\to \infty}\varphi_n(g)=0$ for every $g\in \bigG$, we have that $\lim_{n\to \infty}\varphi_n(a)=0$ for every $a\in A$.

Joining hyperrigidity the property of rigidity at $0$ yields a somewhat stronger notion of hyperrigidity. For example, for a nonunital $C^*$-algebra $A$, these two properties (when combined together) are equivalent to a version of Definition \ref{Dn:HR} in which the faithful $*$-representation $A\subseteq B(H)$ is not assumed to be nondegenerate. In this case, where $A$ is nonunital, the two properties are also equivalent to the hyperrigidity of $\bigG \cup \{1\}$ inside the minimal unitization $A^1$.

We show that while the Cuntz--Krieger family of any directed graph is always rigid at $0$, the smaller generating set $\bigE$ is never rigid at $0$. More precisely, we show that for a row-finite graph, a subset $\bigG$ of the Cuntz--Krieger family is rigid at $0$ if and only if it contains every vertex projection.

One immediate consequence of the latter characterization is that whenever a row-finite graph $G$ has infinitely many vertices (or, equivalently, whenever $\bigO(G)$ is nonunital), $\bigG \cup\{1\}$ is hyperrigid in the unitized $C^*$-algebra $\bigO(G)^1$ if and only if $\bigG$ contains every vertex projection.

\section{Hyperrigidity}

We now describe a few properties of hyperrigidity.
\begin{Pn}\label{Pn:HRprop}
Let $A$ be a $C^*$-algebra and $\bigG$ a generating subset of $A$. Then the following conditions are equivalent:
\begin{enumerate}
\item[(i)] $\bigG$ is hyperrigid in $A$,
\item[(ii)] $\bigG \cup \bigG^*$ is hyperrigid in $A$,
\item[(iii)]${\rm{span}}(\bigG)$ is hyperrigid in $A$.
\end{enumerate} 
If $A$ is unital, then (i)--(iii) are also equivalent to
\begin{enumerate}
\item[(iv)] $\bigG \cup \{1\}$ is hyperrigid in $A$.
\end{enumerate} 
\end{Pn}

\begin{proof}
This follows directly from the definition of hyperrigidity.	
\end{proof}

\begin{Rk}\em
Suppose $A$ is a nonunital $C^*$-algebra generated by a subset $\bigG$. Let $A^1$ denote the minimal unitization of $A$. If $\bigG \cup \{1\}$ is hyperrigid in $A^1$, then $\bigG$ is clearly hyperrigid in $A$. The converse, however, fails; see Proposition \ref{Pn:G1inA1}.
\end{Rk}

The notion of the {\em unique extension property}, originally defined by Arveson, has developed in various settings over the last decades. The common definition \cite[Definition 2.1]{Arveson-note} is for a unital completely positive map defined on an operator system $S$ which generates a unital $C^*$-algebra $A$. The notion has a parallel version for a unital completely contractive map defined on a unital operator algebra, and in \cite[Definition 2.1]{DorOn-Salomon} the latter was generalized to the nonunital case.

In this paper, we will need the following version of the unique extension property which concerns a restriction (to a generating set) of a $*$-representation: the $C^*$-algebra is not assumed to be unital, and the generating set is not assumed to be an operator system or an operator algebra.
\begin{Dn}\label{Dn:UEP}
Let $A$ be a $C^*$-algebra, $\bigG$ a generating subset of $A$ and $\pi:A \to B(H)$ a $*$-representation. We say that $\pi|_{\bigG}$ has the {\em unique extension property} if for every  
completely contractive completely positive map $\rho:A \to B(H)$ we have that
\[
\phi(g)=\pi(g)~~\text{for all } g\in \bigG ~~ \implies ~~ 
\phi(a)=\pi(a)~~\text{for all } a\in A.
\]
\end{Dn}

If $A$ is a unital $C^*$-algebra, $S$ is an operator system generating $A$, and $\pi:A \to B(H)$ is a unital $*$-representation, then 
  $\pi|_{S}$ has the unique extension property if and only if there exists a unique 
{\em unital completely positive} extension of $\pi|_{S}$ to $A$, namely $\pi$ itself.
This shows that if the generating set is an operator system, then our definition for the unique extension property agrees with the one in \cite[Definition 2.1]{Arveson-note}.

Now let $A$ be a nonunital $C^*$-algebra and $\rho:A \to B(H)$ a completely contractive completely positive map.  Let $A^1$ be the minimal unitization of $A$.
The {\em unitization} of $\rho$ is the unital completely positive map $\rho^1:A^1 \to B(H)$ defined by 
\[
\rho^1(a+\lambda \cdot 1)=\rho(a)+\lambda\cdot I_H \quad\text{for all }~ a \in A,~\lambda \in \mathbb C;
\]
see \cite[Proposition 2.2.1]{BrownOzawa} for a proof that $\rho^1$ is indeed a well defined unital completely positive map. If $\rho$ is a $*$-representation, then $\rho^1$ is a $*$-representation as well.

The proof of the following proposition follows immediately from the definition of the unique extension property and the discussion above, and is therefore omitted.
\begin{Pn}\label{Pn:nonunitalUEP}
Let $A$ be a {\em nonunital} $C^*$-algebra, $\bigG$ a generating subset of $A$ and $\pi:A \to B(H)$ a $*$-representation. Then the following conditions are equivalent:
\begin{enumerate}
\item[(i)]  $\pi|_{\bigG}$ has the unique extension property;
\item[(ii)]  $\pi^1|_{\bigG\cup\{1\}}$ has the unique extension property.  
\end{enumerate}
\end{Pn}

In \cite{Arv_HR}, Arveson showed that the unique extension property (in the sense of operator systems and unital $C^*$-algebras) is preserved under direct sums. The proof of the following version is similar and is therefore omitted.
\begin{Pn}\label{Pn:SumOfUep}
Let $A$ be a $C^*$-algebra generated by a subset $\bigG$.
Let $I$ be a set, and for every $i\in I$ let $\pi_i:A\to B(H_i)$ be a $*$-representation such that $\pi_i|_\bigG$ has the unique extension property. Set $\pi:=\oplus_{i\in I} \pi_i$. Then $\pi|_{\bigG}$ has the unique extension property. 
\end{Pn}

The next proposition shows that having the unique extension property for every restriction of a $*$-representation is equivalent to having it only for restrictions of {\em nondegenerate} $*$-representations and for the trivial $*$-representation.

\begin{Pn}\label{Pn:UEPconditions}
Let $A$ be a $C^*$-algebra generated by a subset $\bigG$.
Then the following conditions are equivalent:
\begin{enumerate}[(ii)]
\item[(i)] for every $*$-representation $\pi:A \to B(H)$, $\pi|_{\bigG}$ has the unique extension property;
\item[(ii)] for every {\em nondegenerate} $*$-representation $\pi:A \to B(H)$, $\pi|_{\bigG}$ has the unique extension property, and for the zero map $0:A \to \mathbb C$, $0|_{\bigG}$ has the unique extension property.
\end{enumerate}
If $A$ is unital and $1\in\bigG$, then (i) and (ii) are also equivalent to
\begin{enumerate}
\item[(iii)] for every {\em unital} $*$-representation $\pi:A \to B(H)$, $\pi|_{\bigG}$ has the unique extension property.
\end{enumerate}
If $A$ is nonunital, then (i) and (ii) are also equivalent to
\begin{enumerate}
\item[(iii)$^1$] for every {\em unital} $*$-representation $\pi:A^1 \to B(H)$,  $\pi|_{\bigG\cup \{1\}}$ has the unique extension property.
\end{enumerate}
\end{Pn}

\begin{proof}
(i) $\implies$ (ii) is clear.
For the opposite direction, if $\pi:A \to B(H)$ is any $*$-representation, let $\pi_{nd}:A \to B(K)$ denote its nondegenerate part, so that $\pi=\pi_{nd} \oplus 0^{\oplus \lambda}$ for some cardinality $\lambda$. By Proposition \ref{Pn:SumOfUep} $\pi$ has the unique extension property.

If $A$ is unital and $1\in\bigG$, then for the zero map $0:A\to\mathbb C$, $0|_{\bigG}$ must always have the unique extension property, and nondegenerate $*$-representations are unital, so that (ii) $\iff$ (iii).

If $A$ is nonunital, then the equivalence (i) $\iff$ (iii)$^1$
follows from Proposition \ref{Pn:nonunitalUEP}.
\end{proof}

Let $A$ be a $C^*$-algebra generated by a subset $\bigG$ and let $\pi:A \to B(H)$ be a (perhaps degenerate) $*$-representation. Let $\phi_n:A \to B(H)$ be a sequence of completely contractive completely positive maps satisfying 
\[
\lim_{n \to \infty}\| \phi_n(g) - \pi (g)\|=0~ \text{ for every } g\in\bigG.
\]
A crucial ingredient in the proof of \cite[Theorem 2.1]{Arv_HR} was to define a faithful unital $*$-rep\-re\-sen\-tation $\iota: B(H) \to \ell^\infty(B(H)) / c_0(B(H))$ by
\[
\iota (x) = (x,x,\dots) + c_0(B(H)),
\]
to choose some faithful unital $*$-rep\-re\-sen\-tation $\kappa: \ell^\infty(B(H)) / c_0(B(H)) \to B(K)$, and to show that if $\kappa \circ \iota \circ \pi|_{\bigG}$ has the unique extension property, then 
\[
\lim_{n \to \infty}\| \phi_n(a) - \pi (a)\|=0~ \text{ for every } a\in A.
\]
The original proof of this implication is only for the case where $A$ is unital, $\bigG$ is an operator system, and $\pi$ is unital, but remains the same in our more generalized setting. 

It is worth to note that if $\pi$ is nondegenerate, then $\kappa \circ \iota \circ \pi$ is nondegenerate, and if $\pi$ is the zero map, then $\kappa \circ \iota \circ \pi$ is the zero map. This observation gives rise to the following approximation lemma.

\begin{La} \label{La:Duni}
Let $A$ be a $C^*$-algebra generated by a subset $\bigG$ and let $\pi:A \to B(H)$ be $*$-representation. Let $\phi_n:A \to B(H)$ be a sequence of completely contractive completely positive maps satisfying 
\[
\lim_{n \to \infty}\| \phi_n(g) - \pi (g)\|=0~ \text{ for every } g\in\bigG.
\]
If either
\begin{enumerate}[(i)]
\item every $*$-representation of $A$ has the unique extension property when restricted to $\bigG$; or
\item every nondegenerate $*$-representation of $A$ has the unique extension property when restricted to $\bigG$, and $\pi$ is nondegenerate; or
\item the zero representation of $A$ has the unique extension property when restricted to $\bigG$, and $\pi$ is the zero representation;
\end{enumerate}
then
\[
\lim_{n \to \infty}\| \phi_n(a) - \pi (a)\|=0~ \text{ for every } a\in A.
\]
\end{La}

Proposition \ref{Pn:HRprop} implies that hyperrigidity does not depend on whether $\bigG$ contains the unit of $A$ (if exists) or not. 
There is therefore no loss in assuming that $\bigG$ contains $1$ whenever $A$ is unital.
The following theorem is a generalization of a theorem of Arveson \cite[Theorem 2.1]{Arv_HR}, and Lemma \ref{La:Duni} (ii) is a key ingredient in its proof. The proof essentially consists of minor modifications of the proof of the original theorem and is therefore omitted as well.

\begin{Tm}\label{Tm:ndHR<=>ndUEP}
Let $A$ be a separable $C^*$-algebra and $\bigG$ be a generating subset of $A$.  If $A$ is unital, assume in addition that $1\in\bigG$.
Then $\bigG$ is hyperrigid in $A$ if and only if for every nondegenerate $*$-representation $\pi:A \to B(K)$ on a separable Hilbert space $K$, $\pi|_{\bigG}$ has the unique extension property. 
\end{Tm}

\section{Rigidity at zero}

We now define a closely related notion concerning a certain rigidity property of a generating subset $\bigG$ of a $C^*$-algebra $A$.
\begin{Dn}
\label{Dn:0HR}
Let $A$ be a $C^*$-algebra and $\bigG$ be a generating subset of $A$.
We say that $\bigG$ is {\em rigid at $0$} in $A$ if for every sequence of contractive positive maps $\varphi_n:A \to \mathbb C$, we have
\[
\lim_{n\to\infty} \varphi_n(g)  =0~~\text{for all } g\in \bigG ~~ \implies ~~ 
\lim_{n\to\infty} \varphi_n(a)  =0~~\text{for all } a\in A.
\]
\end{Dn}
\begin{Rk}\label{Rk:1inG=>0HR}\em
If $A$ is a unital $C^*$-algebra and $\bigG$ is a generating subset of $A$ containing the unit, then $\bigG$ must be rigid at $0$ in $A$. To see this, let $\varphi_n:A \to \mathbb C$ be a sequence of contractive positive maps satisfying $\lim_{n\to\infty} \varphi_n(g)  =0$ for all $g\in \bigG$. Then for every $a \in A$
\[
\limsup_{n\to\infty} \|\varphi_n(a)\| \leq \limsup_{n\to\infty} \|\varphi_n(1)\|\|a\| =0.
\]  
This is not the case, however, for nonunital $C^*$-algebras or unital $C^*$-algebras with a generating set that does not contain the unit; see Examples \ref{Ex:Cuntz} and \ref{Ex:Compact}.
\end{Rk}

In the following theorem, we give some equivalent conditions to the property of rigidity at $0$.
\begin{Tm}\label{Tm:0HR}
Let $A$ be a separable $C^*$-algebra generated by a subset $\bigG$. Then the following conditions are equivalent:
\begin{enumerate}[(i)]
\item $\bigG$ is rigid at $0$;
\item for every faithful nondegenerate $*$-representation $A \subseteq B(H)$ and for every sequence of unital completely positive maps $\phi_n:B(H) \to B(H)$, we have
\[
\quad\quad \lim_{n\to\infty}\phi_n(g)= 0~~\text{for all } g\in \bigG ~~ \implies ~~ 
\lim_{n\to\infty}\phi_n(a)= 0~~\text{for all } a\in A;
\]
\item for the zero map $0:A \to \mathbb C$, $0|_{\bigG}$ has the unique extension property;
\item for every separable Hilbert space $K$ and $0:A \to B(K)$, $0|_{\bigG}$ has the unique extension property;
\item there are no states on $A$ vanishing on $\bigG$.
\end{enumerate}
\end{Tm}

\begin{proof}
(v)$\iff$(iii). Clearly, if $0|_{\bigG}$ has the unique extension property, then there are no states on $A$ vanishing on $\bigG$, and converesly, if $0|_{\bigG}$ has a nontrivial positive contractive extension $\varphi$, then $\frac{1}{\|\varphi\|}\varphi$ is a state on $A$ vanishing on $\bigG$.  

(i)$\implies$(iii) is clear.

(ii)$\implies$(v). Assume that (v) does not hold and let $A \subseteq B(H)$ be a faithful nondegenerate $*$-representation. By Arveson's extension theorem we obtain a state $\varphi$ on $B(H)$ vanishing on $\bigG$ but not on $A$. Define $\phi: B(H) \to B(H)$ by $\phi(x):=\varphi(x)  I_H$, then $\phi$ is a unital completely positive map vanishing on $\bigG$ but not on $A$. Setting $\phi_n = \phi$ for every $n \in \mathbb N$ yields a contradiction to (ii).

The assertion (iii)$\implies$(iv) follows by Proposition \ref{Pn:SumOfUep}, while the assertions (iv)$\implies$(ii) and (iv)$\implies$(i) both follow by Lemma \ref{La:Duni} (iii).
\end{proof}

In \cite[Corollary 3.4]{Arv_HR}, Arveson shows that for every $n \in \mathbb N$ the set of $n$ isometries $\bigG_n=\{S_1,\dots, S_n\}$ generating the Cuntz algebra $\bigO_n$ is hyperrigid. 
In \cite[Continuation of Example 2.7]{MuhlySolel98}, Muhly and Solel show that the infinite set of isometries $\bigG_\infty=\{S_1,S_2, \dots\}$ generating the Cuntz algebra $\bigO_\infty$ is not hyperrigid (more precisely, they showed that the norm-closed algebra generated by $\bigG_\infty$ is not hyperrigid).
In the next example we show that for any $n$, finite or not, $\bigG_n$ is not rigid at $0$ in $\bigO_n$.

\begin{Ex}\label{Ex:Cuntz}\em
Let $n\in \mathbb N \cup \{\infty\}$, let $\bigO_n$ be the Cuntz algebra, and let $\bigG_n=\{S_1,S_2,\dots,S_n\}$ be the set of $n$ isometries generating it. We will show that there exists a state on $\bigO_n$ vanishing on $\bigG$. 

For $n=1$ note that $\bigO_1 \cong C(\mathbb T)$, and the state
\[
\rho(f)=\int f(z) d\mu(z),
\]
where $\mu$ is the normalized Lebesgue measure on the unit circle, vanishes on $z$; a precomposition with the above isomorphism yields a state on $\bigO_1$ vanishing on the unitary $S_1$.

Assume now that $n\geq 2$. Let $\bigO_n \subseteq B(H)$ be some faithful nondegenerate $*$-representation, and let $\rho : \bigO_n \to B(H)$ be the completely contractive completely positive map defined by
\[
\theta(a)=S_2^*S_1^* a S_1 S_2, \quad \text{for all }a\in \bigO_n.
\]
Then $\rho|_{\bigG_n}=0$. 
\end{Ex}

\begin{Rk}\em
When $n$ is finite, the study of Davidson and Pitts on Cuntz algebra atomic representations \cite{DavidsonPitts} gives rise to an alternative proof of $\bigG_n$ being not rigid at $0$ in $\bigO_n$. Indeed, suppose that $\pi:\bigO_n \to B(H)$ is any atomic $*$-representation, namely, that there exists an orthonormal basis $\{\xi_j\}$ of $H$, $n$ endomorphisms $\sigma_i:\mathbb N \to \mathbb N$, and scalars $\lambda_{i,j} \in \mathbb T$ such that $\pi(S_i)\xi_j=\lambda_{i,j}\xi_{\sigma_i(j)}$. Let $\xi$ be a wandering vector for the set of noncommutative words in $\pi(\bigG_n)$ (see \cite[Corollary 3.6]{DavidsonPitts}). Then the state defined by
\[
\varphi(a)=\langle \pi (a) \xi, \xi \rangle, \quad \text{for all }a\in \bigO_n
\]   
vanishes on $\bigG_n$.
\end{Rk}

\begin{Pn}\label{Pn:AI=>0HR}
Let $A$ be a $C^*$-algebra and $\bigG$ a generating subset of $A$.
If ${\rm{span}}(\bigG)$ contains an approximate identity for $A$, then it must be rigid at $0$ in $A$. 
\end{Pn}
\begin{proof}
Let $\{e_n\} \subseteq {\rm{span}} (\bigG)$ be an approximate identity for $A$. If $\varphi: A \to \mathbb C$ is a contractive positive linear functional vanishing on $\bigG$, then for every $0\leq a \in A$
\[
\varphi(a)=\varphi(\lim e_n^{\frac{1}{2}} a e_n^{\frac{1}{2}}) \leq \|a\| \lim\varphi(e_n) =0. 
\]
Thus $\varphi=0$, so by Theorem \ref{Tm:0HR} $\bigG$ is rigid at $0$ in $A$.
\end{proof}

\begin{Ex}\label{Ex:CK-family=>0HR}\em
Let $\bigG$ be a directed graph, and let $\bigG:=\{P_v,S_e:v\in V, e\in E\}$ be the universal Cuntz--Krieger family generating $\bigO(G)$. Then ${\rm{span}}(\bigG)$ obviously contains an approximate identity for $\bigO(G)$, namely, the finite sums of the form $\sum_{v\in F} P_v$ where $F$ runs over all finite subsets of $V$. Thus, $\bigG$ is rigid at $0$ in $\bigO(G)$. 
Note that when $G$ is row-finite, this already follows from \cite[Theorem 3.9]{DorOn-Salomon} (together with Theorem \ref{Tm:0HR}).
\end{Ex}

When $A$ is a nonunital $C^*$-algebra generated by a subset $\bigG$, it is natural to ask not only whether $\bigG$ is hyperrigid in $A$, but also whether $\bigG \cup \{1\}$ is hyperrigid in $A^1$; see \cite[Section 6]{Arv_HR}. The answer depends on whether $\bigG$ is rigid at $0$ or not.
\begin{Pn}\label{Pn:G1inA1}
Let $A$ be a nonunital $C^*$-algebra generated by a subset $\bigG$, and let $A^1$ denote its minimal unitization. Then $\bigG \cup\{1\}$ is hyperrigid in $A^1$ if and only if $\bigG$ is rigid at $0$ and hyperrigid in $A$.
\end{Pn}
\begin{proof}
If $\bigG \cup\{1\}$ is hyperrigid in $A^1$, then obviously $\bigG$ is hyperrigid in $A$. As for rigidity at $0$, by Theorem \ref{Tm:ndHR<=>ndUEP} we have that for every unital $*$-representation $\pi:A^1 \to B(H)$, $\pi|_{\bigG}$ has the unique extension property. By Proposition \ref{Pn:UEPconditions}, this
implies that for every (perhaps degenerate) $*$-representation $\pi:A \to B(H)$  --- and, in particular, for the zero representation --- $\pi|_\bigG$ has the unique extension property. By Theorem \ref{Tm:0HR}, $\bigG$ is rigid at $0$ in $A$.

Conversely, assume that $\bigG$ is rigid at $0$ and hyperrigid in $A$. As $\bigG$ is hyperrigid, Theorem \ref{Tm:ndHR<=>ndUEP} implies that for any nondegenerate $*$-representation $\pi:A \to B(H)$, $\pi|_{\bigG}$ has the unique extension property. As $\bigG$ is rigid at $0$, Theorem \ref{Tm:0HR} implies that for the zero representation $0:A \to \mathbb C$, $0|_{\bigG}$ has the unique extension property. By Proposition \ref{Pn:UEPconditions}, for every unital $*$-representation $\pi:A^1 \to B(H)$, $\pi|_{\bigG \cup\{1\}}$ has the unique extension property. Thus, Theorem \ref{Tm:ndHR<=>ndUEP} implies that $\bigG \cup \{1\}$ is hyperrigid in $A^1$.
\end{proof}

In the following theorem we give some equivalent conditions for a set of generators $\bigG$ of a $C^*$-algebra $A$ being both rigid at $0$ and hyperrigid in $A$; we assume, without the loss of generality, that when $A$ is unital $\bigG$ conatins the unit.
\begin{Tm}\label{Tm:UEP<=>HR+0HR}
Let $A$ be a $C^*$-algebra generated by a subset $\bigG$. Then the following conditions are equivalent:
\begin{enumerate}
\item[(i)] $\bigG$ is rigid at $0$ and hyperrigid in $A$;
\item[(ii)] for every $*$-representation $\pi:A \to B(H)$ on a separable Hilbert space $H$, $\pi|_{\bigG}$ has the unique extension property;
\item[(iii)] for every faithful $*$-representation $A \subseteq B(H)$ and every sequence of unital completely positive maps $\phi_n:B(H) \to B(H)$, we have that
\[
\quad\quad\quad\lim_{n\to\infty}\phi_n(g)= g~~\text{for all } g\in \bigG ~~ \implies ~~ 
\lim_{n\to\infty}\phi_n(a)= a~~\text{for all } a\in A
\]
where all convergences are in norm.
\end{enumerate}
If $A$ is unital , then (i)--(iii) are also equivalent to
\begin{enumerate}
\item[(iv)] $\bigG$ is hyperrigid in $A$.
\end{enumerate}
If $A$ is nonunital, then (i)--(iii) are also equivalent to
\begin{enumerate}
\item[(iv)$^1$] $\bigG \cup \{1\}$ is hyperrigid in $A^1$.
\end{enumerate}
\end{Tm} 

\begin{proof}
If $A$ is unital and $1\in \bigG$, then $\bigG$ is always rigid at $0$ (see Remark \ref{Rk:1inG=>0HR}), so in this case we have (i) $\iff$ (iv). 
If $A$ is nonunital, then by Proposition \ref{Pn:G1inA1} we have (i) $\iff$ (iv)$^1$.

(i) $\implies$ (ii). Assume that $\bigG$ is both rigid at $0$ and hyperrigid in $A$. As $\bigG$ is rigid at $0$ by Theorem \ref{Tm:0HR} we have that for the zero representation $0:A \to \mathbb C$, $0|_{\bigG}$ has the unique extension property. As $\bigG$ is hyperrigid, by Theorem \ref{Tm:ndHR<=>ndUEP} we have that for any nondegenerate $*$-representation $\pi:A \to B(H)$, $\pi|_{\bigG}$ has the unique extension property. By the equivalence (i) $\iff$ (ii) of Proposition \ref{Pn:UEPconditions}, we are done.

(ii) $\implies$ (iii). Follows by Lemma \ref{La:Duni} (i).

(iii) $\implies$ (i). We clearly have that condition (iii) implies hyperrigidity, so we only need to show it implies rigidity at $0$. 
We will do so by obtaining the equivalent condition (ii) from Theorem \ref{Tm:0HR}. 
Let $\iota: A  \hookrightarrow B(H)$ be a faithful nondegenerate $*$-representation and let $\phi_n:B(H) \to B(H)$ be a sequence of unital completely positive maps satisfying $\lim_{n \to \infty} \phi_n(g)=0$ for every $g \in \bigG$. Consider the faithful (degenerate) $*$-representation $\iota \oplus 0 : A \hookrightarrow B(H\oplus H)$ and identify $A$ with $\iota \oplus 0 (A)$. For every $n\in \mathbb N$, let $\psi_n$ denote the extension to $B(H\oplus H)$ of ${\rm{id}}_{B(H)}\oplus \phi_n$. Then $\lim_{n \to \infty} \psi_n(g)=g$  for every $g \in \bigG$, so $\lim_{n\to \infty} \psi_n(a)=a$ for every $a \in A$. Thus, 
$\lim_{n \to \infty} \phi_n(a)=0$ for every $a\in A$.
\end{proof}


\section{Hyperrigidity of the edge set}
Let $G=(V,E,s,r)$ be a directed graph. 
Recall that the Cuntz--Krieger algebra $\bigO(G)$ is the universal $C^*$-algebra generated by a Cuntz--Krieger $G$-family $(\mathcal V,\mathcal E)$, where $\mathcal V:= \{P_v : v \in V\}$ and $\mathcal E:= \{S_e : e \in E\}$ satisfying the relations (I), (TCK) and (CK) described in the Introduction. We sometimes call $\mathcal V$ and $\mathcal E$ the vertex set and the edge set of $\bigO(G)$, respectively. As was mentioned in the Introduction (see \cite[Theorems 3.5 and 3.9]{DorOn-Salomon}), it is known that $\mathcal V \cup \mathcal E$ is hyperrigid in $\bigO(G)$ if and only if $G$ is row-finite, so if $G$ is non-row-finite, $\bigE$ can never be hyperrigid in $\bigO(G)$.

If $G$ is row-finite, then $\bigE$ generates $\bigO(G)$ if and only if $G$ contains no isolated vertices (namely, there are no vertices that emit and receive no edges). In this section, we show that in this case $\bigE$ is hyperrigid in $\bigO(G)$. We start with a finite graph and then continue, by taking inductive limits, to any row-finite graph.

\subsection{Finite graphs}
A directed graph $G=(V,E,s,r)$ is called {\em finite} if both $V$ and $E$ are finite.
To show that for finite graphs with no isolated vertices $\bigE$ is hyperrigid in $\bigO(G)$, we will use the machinery of maximal dilations of unital completely positive maps on operator systems.

Let $S$ be an operator system. A unital completely positive map $\phi:S \to B(H)$ is said to be {\em maximal} if whenever $\psi: S \to B(K)$ is a unital completely positive map dilating $\phi$ --- that is, $K \supseteq H$ and $\phi=P_H \psi(~\cdot~) |_{H}$ ---  then $\psi=\phi \oplus \rho$ for some unital completely positive map $\rho$. It is known that any completely positive map $\phi:S \to B(H)$
dilates to a maximal map; see \cite[Theorem 1.3]{Arveson-note} (this was originally proved in \cite{DritschelMcCullough}, but in terms of operator algebras rather than operator systems). 

The notions of maximality and the unique extension property are strongly related. 
Suppose $A$ is a unital $C^*$-algebra generated by an operator system $S$, then a unital completely positive map $\phi:S \to B(H)$ is maximal if and only if it extends uniquely to a unital completely positive map on $A$ and additionally this extension is $*$-rep\-re\-sen\-tation. In particular, if $\pi:A \to B(H)$ is a unital $*$-representation, then $\pi|_{S}$ is maximal if and only if it has the unique extension property; see \cite[Proposition 2.2]{Arveson-note} (the original proof appeared in \cite{MuhlySolel_boundary}, but again in terms of operator algebras rather than operator systems). 

\begin{Tm}\label{Tm:HR_edges}
Let $G=(V,E,s,r)$ be a finite graph with no isolated vertices. Then $\bigE := \{S_e: e\in E\}$ is hyperrigid in $\bigO(G)$.
\end{Tm}
\begin{proof}
Let $\pi: \bigO(G) \to B(H)$ be a unital $*$-representation. By Proposition \ref{Pn:HRprop} and Theorem \ref{Tm:ndHR<=>ndUEP} we need to show that $\pi|_{\bigE + \bigE^* + \mathbb C 1}$ has the unique extension property, or equivalently, that it is maximal. Let $\widetilde \rho$ be a maximal dilation of $\pi|_{\bigE + \bigE^* + \mathbb C 1}$ to a Hilbert space $K \supseteq H$, and let $\rho:\bigO(G) \to B(K)$ be its extension to a $*$-representation. 
Denote
\[
\rho(S_e)=
\begin{bmatrix}
\pi(S_e) & X_e\\
Y_e & Z_e
\end{bmatrix} \quad \text{for all } e\in E.
\]
Let $W$ be a subset of $V$ containing all sources and no sinks.
For every $v \in W$ choose a representative $e_v \in s^{-1}(v)$.
Then
\[
1
=\sum_{v\in W} P_v + \sum_{v\in V\setminus W} P_v
=\sum_{v\in W} S_{e_v}^*S_{e_v} + \sum_{v\in V\setminus W} \sum_{r(e)=v} S_e S_e^*. 
\]

Thus,
\[
\begin{split}
1
&=\sum_{v\in W} \rho(S_{e_v})^*\rho(S_{e_v}) + \sum_{v\in V\setminus W} \sum_{r(e)=v} \rho(S_e) \rho(S_e)^*\\
&=
\sum_{v\in W}
\begin{bmatrix}
\pi(S_{e_v})^* \pi(S_{e_v}) +Y_{e_v}^* Y_{e_v} & *\\
* & *
\end{bmatrix} \\
&\quad\quad+
\sum_{v\in V\setminus W}
\sum_{r(e)=v}
\begin{bmatrix}
\pi(S_e) \pi(S_e)^* +X_e X_e^* & *\\
* & *
\end{bmatrix},
\end{split}
\]
so $\sum_{v\in W} Y_{e_v}^*Y_{e_v} + \sum_{v\in V\setminus W} \sum_{r(e)=v} X_e X_e^*=0$.
We therefore conclude that $X_e=0$ for all $e\in r^{-1}(V \setminus W)$, and $Y_{e_v}=0$ for all $v \in W$. Since the representatives $\{e_v\}_{v\in W}$ were chosen arbitrarily, we have $Y_{e} =0$ for all $e\in s^{-1}(W)$.

Finally, note that choosing $W$ to be the set of all sources in the graph, we obtain that $r^{-1}(V\setminus W)=E$, so $X_e=0$ for {\em all} $e\in E$; similarly choosing $V\setminus W$ to be the set of all sinks in the graph, we obtain that $s^{-1}(W)=E$, so $Y_e=0$ for {\em all} $e\in E$. Thus, $\pi|_{\bigE + \bigE^* + \mathbb C 1}$ is maximal.
\end{proof}

\subsection{Row-finite graphs}
We will now show that for a row-finite directed graph $G$, whenever $\bigE$ generates $\bigO(G)$, it must be hyperrigid. For this, we first recall that a row-finite graph $G$ is an inductive limit (or equivalently, a {\em direct union}) of a sequence of finite subgraphs, and consequently that $\bigO(G)$ is the inductive limit of the corresponding finite graph Cuntz--Krieger algebras. These results are considered as folklore, but we provide the details for completeness.

Let $G_1=(V_1,E_1,s_1,r_1)$ and $G_2=(V_2,E_2,s_2,r_2)$ be two directed graphs. Then $G_1$ is said to be a {\em subgraph} of $G_2$ if $V_1 \subseteq V_2$, $E_1 \subseteq E_2$, $s_1=s_2|_{E_1}$ and $r_1=r_2|_{E_1}$. If furthermore whenever $v \in V_1$ is not a source or an infinite receiver in $G_1$ we have $r_1^{-1}(v)=r_2^{-1}(v)$, then $G_1$ is said to be a {\em CK subgraph} of $G_2$.

Let $\{ G_\alpha=(V_\alpha,E_\alpha,s_\alpha,r_\alpha) \}_{\alpha \in A }$ be a family of CK subgraphs of $G=(V,E,s,r)$ {\em directed under CK inclusion}, i.e., $(A,\leq)$ is a directed set of indices, and whenever $\alpha \leq \beta$ we have that $G_\alpha$ is a CK subgraph of $G_\beta$. Then $G = (V,E,r,s)$ is said to be the {\em direct union} of the family $\{G_\alpha\}_{\alpha \in A}$ if $V=\bigcup_\alpha V_\alpha$, $E=\bigcup_\alpha E_\alpha$, and for every $\alpha \in A$, $e\in E_\alpha$ we have $r(e) := r_{\alpha}(e)$, $s(e) := s_{\alpha}(e)$.

\begin{Pn}\label{Pn:row-finite<=>indlimfin}
Let $G=(V,E,s,r)$ be a directed graph. Then $G$ is row-finite if and only if it is a direct union of finite CK subgraphs. If $G$ is a row-finite graph containing no isolated vertices, then we can choose the finite CK subgraph to contain no isolated vertices as well.
\end{Pn}
\begin{proof}
Suppose $G$ is row-finite. For every $X$ finite subset of $V$, let $E_X$ consists of all edges with $r(e)\in X$ (by the row-finiteness assumption $E_X$ is finite) and let $V_X:=X\cup s(E_X)$ (so $V_X$ is finite as well).
Let $s_X:=s|_{E_X}$, $r_X:=r|_{E_X}$, and $G_X:=(V_X, E_X, s_X, r_X)$.
Clearly, $G_X$ is a subgraph of $G$. Since each $v\in V_X \setminus X$ is a source in $G_X$, and each $v\in X$ has $r_X^{-1}(v)=r^{-1}(v)$, we have that $G_X$ is a CK subgraph of $G$. A similar argument shows that if $X \subseteq Y$, then $G_X$ is a CK subgraph of $G_Y$. Thus, $G$ is the inductive union of the family $\{G_X\}$ indexed over all finite subsets $X\subset V$,  which is directed under CK inclusion.

If $G$ is not row-finite, there is an infinite receiver $v\in V$, and any finite CK subgraph of $G$ cannot contain all the edges $r^{-1}(v)$. Therefore, $v$ is a source for \emph{every} finite CK subgraph of $G$ containing it. Indeed, if $G_1=(V_1,E_1,s_1,r_1)$ is a finite CK subgraph of $G$ containing $v$, and the latter is not a source for $G_1$, then $r_1^{-1}(v) = r^{-1}(v)$, which is impossible as $r^{-1}(v)$ is infinite. Thus, $v$ must be a source in any union of finite CK subgraphs of $G$, so $G$ is not the union of finite CK subgraphs. 

Suppose now that $G$ contains no isolated vertices, and let $S \subseteq V$ be the subset of all sources in $G$. Since $G$ contains no isolated vertices, we can choose for each $v\in S$ a representative 
$e_v \in s^{-1}(v)$. For every $X$ finite subset of $V$, let $X'=X\cup\{r(e_v) : v\in S\cap X\}$ (so $X'$ is finite as well). Let $G_{X'}=(V_{X'},E_{X'},s_{X'},r_{X'})$ be the CK subgraph of $G$ described in the first paragraph of the proof (constructed with respect to $X'$). 
Then clearly, $G$ is still the inductive union of the family $\{G_{X'}\}$ where $X$ is a finite subset of $V$. 
Assume that $V_{X'}$ contain an isolated vertex $v$ of $G_{X'}$. As $v$ is isolated it is not in $s(E_{X'})$, so $v\in X'$. But $E_{X'}=r^{-1}(X')$, so that $v$ must be a source (of $G$). But then $v\in S\cap X'=S\cap X$, so that $r(e_v) \in X'$. Thus, $e_v \in E_{X'}$ and $s(e_v)=v$, which is a contradiction.
\end{proof} 

\begin{Pn}\label{Pn:CK_subgraph}
If $G_1=(V_1,E_1,s_1,r_1)$ is a CK subgraph of a directed graph $G_2=(V_2,E_2,s_2,r_2)$, then there is an embedding of $C^*$-algebras $\bigO(G_1) \hookrightarrow \bigO(G_2)$ mapping generators to generators.
\end{Pn}
\begin{proof}
As $G_1$ is a subgraph of $G_2$ we have $V_1 \subseteq V_2$ and $E_1 \subseteq E_2$. We show that the Cuntz--Krieger family $\{P_v,S_e:v\in V_1, e\in E_1\} \subseteq \bigO(G_2)$ of $G_2$ is also a Cuntz--Krieger family for $G_1$. Condition (I) is clearly satisfied. As for (CK), let $v \in V_1$. If $v$ is either a source in $G_1$ or an infinite receiver, then there is nothing to check. Otherwise, since $G_1$ is a {\em CK} subgraph in $G_2$, we have $r_1^{-1}(v)=r_2^{-1}(v)$, so $\sum_{e \in r_1^{-1}(v)}S_eS_e^*= \sum_{e \in r_2^{-1}(v)}S_eS_e^* = P_v$, and (CK) is satisfied as well. Thus, by the universality of $\bigO(G_1)$, there exists a $*$-homomorphism $\varphi:\bigO(G_1) \to \bigO(G_2)$, and by the gauge-invariant uniqueness theorem, $\varphi$ must be injective.
\end{proof}

We now obtain the following corollary.

\begin{Cy}\label{Cy:CK_indlim}
If $G$ is a row-finite graph, then there is a directed family $\{G_\alpha:\alpha \in A\}$ of finite CK subgraphs of $G$ such that the union of $\bigO(G_\alpha)$ is dense in $\bigO(G)$.
If $G$ is a row-finite graph containing no isolated vertices, then we can choose the finite CK subgraph $G_\alpha$ to contain no isolated vertices as well.
\end{Cy}

We now show that hyperrigidity is preserved by taking inductive limits with injective connecting maps.

\begin{Pn}\label{Pn:HR_indlim}
Let $A$ be a $C^*$-algebra and $\bigG$ a generating subset of $A$. Suppose that there is a collection $\{A_\alpha\}$ of $C^*$-subalgebras of $A$ directed under inclusion and with a dense union, and let $\bigG_\alpha :=\bigG \cap A_\alpha$. If for every $\alpha$ we have that $\bigG_\alpha$ is hyperrigid in $A_\alpha$, then $\bigG$ is hyperrigid in $A$.
\end{Pn}

\begin{proof}
Let $\pi:A \to B(H)$ be nondegenerate $*$-representation.
For every $\alpha$, let $\pi_\alpha: A_\alpha \to B(H_\alpha)$ be the nondegenerate part of $\pi|_{A_\alpha}:A_\alpha\to B(H)$, that is, $H_\alpha= \overline{\pi(A_{\alpha})H}$ and $\pi_\alpha= P_{H_\alpha}\pi|_{A_{\alpha}}(\cdot)|_{H_\alpha}$.
By assumption, $\pi_\alpha|_{\bigG_\alpha}$ has the unique extension property. 

Let $\varphi : A \to B(H)$ be a completely contractive completely positive extension of $\pi|_{\bigG}$, and consider
$\varphi|_{A_\alpha}:A_\alpha\to B(H)$. Then $x \mapsto P_{H_\alpha}\varphi|_{A_\alpha}(x) |_{H_\alpha}$ is a completely contractive completely positive map $A_\alpha \to B(H_\alpha)$, and for every $g\in \bigG_\alpha$ we have 
\[
P_{H_\alpha}\varphi|_{A_\alpha}(g) |_{H_\alpha}=P_{H_\alpha}\pi(g) |_{H_\alpha}=\pi_\alpha(g).
\] 
Thus, $\pi(a)=P_{H_\alpha}\varphi(a) |_{H_\alpha}$ for every $a \in A_\alpha$.
Now for every $a\in A_{\alpha}$ and $\xi \in H_{\beta}$ choose $\gamma\geq \alpha,\beta$. Then,
\[
\pi(a)\xi=P_{H_\gamma} \pi(a) |_{H_\gamma}\xi=\pi_\gamma(a)\xi=P_{H_\gamma}\varphi(a) P_{H_\gamma}\xi=P_{H_\gamma}\varphi(a)\xi.
\]
Since this holds for any $\gamma \geq \alpha,\beta$ and since $\bigcup H_{\gamma}$ is dense in $H$, we have that
$
\pi(a)\xi=\varphi(a)\xi.
$
As this is true for every $\xi \in H_{\beta}$ and $a\in A_{\alpha}$ and for every $\alpha$ and $\beta$ we conclude that $\pi=\varphi$.
\end{proof}

Proposition \ref{Pn:HR_indlim} together with Theorem \ref{Tm:HR_edges} and Corollary \ref{Cy:CK_indlim} implies that for row-finite graphs, whenever $\bigE$ generates $\bigO(G)$, then it must be hyperrigid. 
\begin{Tm}\label{Tm:RF=>ndHR}
Let $G$ be a row-finite graph with no isolated vertices. Then $\bigE$ is hyperrigid in $\bigO(G)$.
\end{Tm}

Note that while for row-finite graphs an equivalent condition for $\bigE$ to generate $\bigO(G)$ is that the graph contains no isolated vertices, in the not necessarily row-finite case one needs to add to the latter condition the requirement that every infinite receiver is also an emitter. 

\begin{Cy}\label{Cy:main}
Let $G$ be a directed graph with no isolated vertices and in which every infinite receiver is an emitter. Then $\bigE$ is hyperrigid in $\bigO(G)$ if and only if $\bigG$ is row-finite.
\end{Cy}

\begin{proof}
If $G$ is not row-finite, then by (the proof of) \cite[Theorem 3.9]{DorOn-Salomon}, even the larger set $\{P_v,S_e:v\in V,e\in E\}$ is not hyperrigid in $\bigO(G)$, so that $\bigE$ itself obviously cannot be  hyperrigid.
If $G$ is row-finite, then Theorem \ref{Tm:RF=>ndHR} implies that $\bigE$ is hyperrigid in $\bigO(G)$.
\end{proof}


Let $A$ be an operator algebra. A {\em $C^*$-cover} of $A$ is a pair $(B,\iota)$ consisting of a $C^*$-algebra $B$ and a completely isometric homomorphism $\iota:A \to B$ such that $\iota(A)$ generates $B$. 
Arveson defined the notion of the $C^*$-envelope. Let $A$ be a unital operator algebra. A $C^*$-cover $(C^*_e(A), \kappa)$ is called a {\em $C^*$-envelope} of $A$ if for every other $C^*$-cover $(B,\iota)$ there exists a (necessarily unique and onto) $*$-homomorphism $\pi:B \to C^*_e(A)$ such that $\pi\circ \iota = \kappa$. In this case, 
$(C^*_e(A), \kappa)$ must be unique up to $*$-isomorphism. Hamana proved that every unital operator algebra admits a $C^*$-envelope \cite[Theorem 4.4]{Hamana}. 

If an operator algebra $A$ is nonunital, then a theorem of Meyer \cite[Corollary 2.1.15]{BlecherLeMerdy} shows that it admits a unique minimal unitization. More precisely, if $(\iota, B)$ is a $C^*$-cover for the operator algebra $A$, and $B \subseteq B(H)$ is some faithful nondegenerate $*$-representation of $B$, then the operator algebraic structure of $A^1 \cong \iota(A) + \mathbb{C} I_H$ is independent of the $C^*$-cover $B$ and of the faithful nondegenerate $*$-representation $B \subseteq B(H)$.

In this case, the $C^*$-envelope of $A$ is defined as follows: let $(C^*_e(A^1),\kappa)$ be the $C^*$-envelope of $A^1$, then the $C^*$-envelope of $A$ is $(C^*(\kappa(A)), \kappa)$ where $C^*(\kappa(A))$ is the $C^*$-algebra generated by $\kappa(A)$ in $C^*_e(A^1)$. It satisfies a universal property generalizing the one for $C^*$-envelopes of unital operator algebras: for any other $C^*$-cover $(B,\iota)$ of $A$ there exists a (necessarily unique and onto) $*$-homomorphism $\pi:B \to C^*_e(A)$ such that $\pi\circ \iota = \kappa$. 

The following proposition and corollary are known for unital operator algebras; see \cite[Theorem 4.1]{DritschelMcCullough}, \cite[Section 2.2]{Arveson_C-starI}.
\begin{Pn}\label{Pn:UEP=>C-star-envelope}
Let $B$ be a $C^*$-algebra generated by an operator algebra $A$. Let $\iota:B\to B(H)$ be some nondegenerate faithful $*$-representation. If $\iota|_A$ has the unique extension property, then $B$ is the $C^*$-envelope of $A$.
\end{Pn} 

\begin{proof}
If $A$ is unital, then this follows from the proof of \cite[Theorem 4.1]{DritschelMcCullough}.
If $B$ is nonunital, then by Proposition \ref{Pn:nonunitalUEP}, for the unital faithful $*$-representation $\iota^1:B^1 \to B(H)$, $\iota^1|_{A^1}$ has the unique extension property. Thus, by the unital case
$C^*_e(A^1)=B^1$, so that $C^*_e(A)$ is the $C^*$-algebra generated by $A$ in $B^1$, namely, $B$.

Finally, if $B$ is unital and $A$ is not, then as $\iota|_A$ has the unique extension property, $\iota|_{A^1}$ must have the unique extension property as well. Thus, $C^*_e(A^1)=B$, so that $C^*_e(A)$ is the $C^*$-algebra generated by $A$ in $B$, namely, $B$ itself.
\end{proof}

\begin{Cy}
Let $B$ be a $C^*$-algebra generated by an operator algebra $A$. If $A$ is hyperrigid in $B$, then $C^*_e(A)=B$.
\end{Cy}

\begin{proof}
If $B$ is nonunital or $A$ is unital, then this follows from Proposition \ref{Pn:UEP=>C-star-envelope} together with Theorem \ref{Tm:ndHR<=>ndUEP}.
Otherwise, $B$ is unital and $A$ is not. If $A$ is hyperrigid in $B$, then $A^1$ is hyperrigid in $B$. Thus, by the unital case, $C^*_e(A^1)=B$, so that $C^*_e(A)$ is the $C^*$-algebra generated by $A$ in $B$, namely, $B$ itself.
\end{proof}

Katsoulis and Kribs showed in  \cite[Theorem 2.5]{KatsoulisKribs_Graphs} that $\bigO(G)$ is the $C^*$-envelope of the norm-closed algebra generated by the Cuntz--Krieger family --- namely, the tensor-algebra $\bigT_+(G)=\overline{{\rm{alg}}}\{P_v,S_e:v\in V, e\in E\}$ (in fact, the tensor-algebra is usually defined as the norm-closed algebra generated by the universal {\em Toeplitz}--Cuntz--Krieger, but as $\bigO(G)$ is a $C^*$-cover of the latter algebra, the two definitions are equivalent).  Our hyperrigidity result shows that for a row-finite graph $G$ the smaller operator algebra $A(\bigE):=\overline{{\rm{alg}}}\{S_e:e\in E\}$ is hyperrigid in $\bigO(G)$, and therefore, the latter must be the $C^*$-envelope of $A(\bigE)$.
\begin{Cy}\label{Cy:c-star-envelope-of-A}
Let $G$ be a row-finite graph containing no isolated vertices, and let  $A(\bigE):=\overline{{\rm{alg}}}\{S_e:e\in E\}$ be the operator algebra generated by $\bigE$ inside $\bigO(G)$. Then,
$
C^*_{e}(A(\bigE))=\bigO(G).
$
\end{Cy}

\section{Rigidity at zero of the edge set}

Let $G=(V,E,s,r)$ be a directed graph.
As was mentioned in the Introduction and in Example \ref{Ex:CK-family=>0HR}, the Cuntz--Krieger $G$-family $\bigV \cup \bigE$ is hyperrigid if and only if $G$ is row-finite, and is rigid at $0$ in any case.
In the last section we showed that whenever $\bigE$ generates $\bigO(G)$, then it is hyperrigid in $\bigO(G)$ if and only if $G$ is row-finite. In this section we show that $\bigE$ is never rigid at $0$ in $\bigO(G)$. More precisely, we show that a generating subset $\bigG \subseteq \bigV\cup\bigE$ is rigid at $0$ if and only if $\bigV \subseteq \bigG$.

\begin{Tm}\label{Tm:Not0HR}
Let $G=(V,E,s,r)$ be a directed graph and let $\bigG$ be a subset of the Cuntz--Krieger family $\{P_v,S_e:v\in V, e\in E\}$. If there exists a vertex $v\in V$ such that $P_v \not \in \bigG$, then $\bigG$ is not rigid at $0$ in $\bigO(G)$.
\end{Tm}
\begin{proof}
Let $v$ be a vertex in $V$ such that $P_v \not \in \bigG$ and let $\bigO(G) \subseteq B(H)$ be a nondegenerate faithful $*$-representation. By Theorem \ref{Tm:0HR} it suffices to show that for the zero map $0:\bigO(G) \to B(H)$, $0|_{\bigG}$ does not have the unique extension property.

If there are no edges $e\in E$ with $r(e)=s(e)=v$ --- namely, if there are no loops at $v$ --- then the nontrivial completely contractive completely positive map $\rho:\bigO(G) \to B(H)$ defined by $\rho(a)=P_vaP_v$ vanishes on $\bigG$, and we are done. On the other hand, if there are two (or more) loops, $e$ and $f$, at $v$, then the nontrivial completely contractive completely positive map $\rho:\bigO(G) \to B(H)$ defined by $\rho(a)=S_f^*S_e^*P_vaP_vS_eS_f$ vanishes on $\bigG$. Thus, we may assume that there is exactly one loop at $v$, which we denote by $e$.

If there exists an edge $f\in E\setminus \{e\}$ with $r(f)=v$, then the nontrivial completely contractive completely positive map $\rho:\bigO(G) \to B(H)$ defined by $\rho(a)=S_f^*P_vaP_vS_f$ vanishes on $\bigG$.

Thus, the only case we need to consider is when there is only one edge going into $v$, namely, $e$. In this case $S_{e}^*S_{e}=P_v=S_eS_e^*$, so there exists a $*$-homomorphism $\pi: C(\mathbb T)\to C^*(S_e)$ mapping 
$z$ to $S_e$. As $\bigO(G)$ admits a natural gauge action and $\pi$ must be equivariant with respect to this gauge action, the gauge invariant uniqueness theorem implies that $\pi$ is faithful.
Thus, by Example \ref{Ex:Cuntz}, there exists a state on $C^*(S_e)$ vanishing on $\{S_e\}$. 
By Arveson's extension theorem it extends to a state $\varphi$ on $P_v \bigO(G) P_v$. 
The nontrivial completely contractive completely positive map $\rho : \bigO(G) \to B(H)$ defined by
$
\rho(a):=\varphi ( P_v a P_v)  I_H
$
must vanish on $\bigG$.
\end{proof}

\begin{Rk}\em
As the referee pointed out, one can prove Theorem \ref{Tm:Not0HR} alternatively using an integration over the gauge action. More precisely, if $P_v \not \in \bigG$ and $\gamma$ denotes the natural gauge action, then an integration of $P_v \gamma_z(\cdot )P_v$ over the unit circle gives rise to a nontrivial completely contractive completely positive extension of $0|_{\bigG}$. 
\end{Rk}

Recall that $\bigO(G)$ is unital if and only if $G$ contains only finitely many vertices. In case $\bigO(G)$ is nonunital, then due to Proposition \ref{Pn:G1inA1}, whenever $\bigG$ is not rigid at $0$, $\bigG\cup \{1\}$ is not hyperrigid in $\bigO(G)^1$.  
\begin{Cy}
Let $G=(V,E,s,r)$ be a directed graph with infinitely many vertices and let $\bigG$ be any subset of the Cuntz--Krieger family. If there exists a vertex $v\in V$ such that $P_v\not\in \bigG$, then $\bigG \cup \{1\}$ is not hyperrigid in $\bigO(G)^1$.
\end{Cy}

\begin{Cy}\label{Cy:Graphs0HR}
Let $G$ be a directed graph with no isolated vertices and in which every infinite receiver is an emitter (or equivalently assume that $\bigO(G)$ is generated by $\bigE$). Then $\bigE$ is not rigid at $0$ in $\bigO(G)$.
\end{Cy}

\begin{Ex}\label{Ex:Compact}\em
Let $n\in \mathbb N \cup \{\infty\}$. Define a directed graph by setting $V=\{i: 1\leq i <n+1 \}$ and $E=\{e_i : 2\leq i < n+1 \}$ where for each $i \in V$ we have $s(e_i)=i$ and $r(e_i)=i-1$.
Set
\[
\begin{split}
E_{i,j}&:=S_{e_{i+1}} S_{e_{i+2}}\cdots S_{e_{j}}, \quad \text{for every } i<j;\\
E_{i,j}&:=E_{j,i}^*, \quad \text{for every } i>j; ~\text{and}\\
E_{i,i}&:=E_{i,1}E_{1,i}=P_i.
\end{split}
\]
A simple computation shows that the $E_{i,j}$'s are matrix units, so that $\bigO(G)=C^*(\bigE)$ can be identified with the $C^*$-algebra $\bigK(H)$ of compact operators on a separable Hilbert space $H$ (in case $n$ is finite, then we obtain $M_n(\mathbb C)$). By Corollaries \ref{Cy:main} and  \ref{Cy:Graphs0HR}, $\bigE=\{E_{i-1,i}:2\leq i < n+1\}$ is hyperrigid but not rigid at $0$ in $\bigO(G)$. 
In addition, note that if $n=\infty$, by Proposition \ref{Pn:G1inA1}, $\bigE \cup \{1\}$ is not hyperrigid in $K(H)^1$.
\end{Ex}

To conclude this section we give the following observation for the case where the edge set $\bigE$ generates $\bigO(G)$. As we showed, in this case, $\bigE$ is hyperrigid if and only if $G$ is row-finite, and $\bigE$ is never rigid at $0$.
The set $\bigE\cup \bigE^*\bigE$  is also hyperrigid if and only if $G$ is row-finite and is always rigid at $0$. To complete the picture, we show that for any graph $G$, row-finite or not, $\bigE\cup \bigE^*\bigE \cup \bigE \bigE^*$ is hyperrigid (and rigid at $0$).
More precisely, we have the following theorem.
\begin{Tm}
Let $G=(V,E,s,r)$ be a directed graph. Then the set $\{P_v, S_e, S_eS_e^*: v\in V, e\in E\}$ is hyperrigid and rigid at $0$ in $\bigO(G)$.
\end{Tm}
\begin{proof}
First assume that $V$ is infinite.
By Propositions \ref{Pn:G1inA1} and \ref{Pn:HRprop} we must show that the operator system $S={\rm{span}}\{1,P_v, S_e, S_eS_e^*: v\in V, e\in E\}$ is hyperrigid in $\bigO(G)^1$. To this end, we use Theorem \ref{Tm:ndHR<=>ndUEP}.
Let $\pi: \bigO(G)^1 \to B(H)$ be a unital $*$-representation. We will show that $\pi|_{S}$ has the unique extension property, or equivalently, that it is maximal. Let $\widetilde \rho: S \to B(K)$ be a maximal dilation of $\pi|_{S}$, and let $\rho:\bigO(G)^1 \to B(K)$ be its extension to a $*$-representation. 
Denote
\[
\rho(S_e)=
\begin{bmatrix}
\pi(S_e) & X_e\\
Y_e & Z_e
\end{bmatrix} \quad \forall e\in E,
\]
\[
\rho(P_v)=
\begin{bmatrix}
\pi(P_v) & X_v\\
Y_v & Z_v
\end{bmatrix} \quad \forall v\in V,
\]
and
\[
\rho(S_eS_e^*)=
\begin{bmatrix}
\pi(S_eS_e^*) & X'_e\\
Y'_e & Z'_e
\end{bmatrix} \quad \forall e\in E.
\]
Let $P:K\to H$ denote the orthogonal projection of $K$ onto $H$, then for every $v\in V$
\[
\begin{split}
P \rho(P_v)^*(1-P) \rho(P_v) P &= P \rho(P_v) P - P \rho(P_v) P \rho(P_v) P\\
&=\pi(P_v) - \pi(P_v)\pi(P_v) = 0,
\end{split}
\]
and the $C^*$-identity implies $Y_v = (1-P) \rho(P_v) P = 0$ for every $v\in V$. As $\rho(P_v)$ is self-adjoint, we have $X_v = 0$ for every $v\in V$ as well. 
A similar argument shows that $X'_e=Y'_e=0$ for every $e \in E$.
Now, for all $e\in E$, we have that $S_e^*S_e=P_{s(e)}$ so that
\[
\begin{bmatrix}
\pi(P_{s(e)}) & 0\\
0 & *
\end{bmatrix}
=\rho(S_e^*S_e)=\rho(S_e)^*\rho(S_e)=
\begin{bmatrix}
\pi(S_e)^*\pi(S_e)+Y_e^*Y_e & *\\
* & *
\end{bmatrix},
\]
which implies $Y_e=0$ for all $e\in E$.

Finally, for all $e\in E$ we have 
\[
\begin{bmatrix}
\pi(S_eS_e^*) & 0\\
0 & *
\end{bmatrix}
=\rho(S_eS_e^*)=\rho(S_e)\rho(S_e)^*=
\begin{bmatrix}
\pi(S_e)\pi(S_e)^*+X_eX_e^* & *\\
* & *
\end{bmatrix},
\]
and we have that $X_e=0$ for all $e \in E$.
Thus $\pi$ is maximal.

If $V$ is finite, then a similar argument shows that the operator system $S={\rm{span}}\{P_v, S_e, S_eS_e^*: v\in V, e\in E\}$ is hyperrigid and rigid at $0$ in $\bigO(G)$.
\end{proof}

To summarize our results, we state the following theorem.
\begin{Tm}
Let $G$ be a directed graph with no isolated vertices and in which every infinite receiver is an emitter (or equivalently assume that $\bigO(G)$ is generated by $\bigE$). Then inside $\bigO(G)$
\begin{enumerate}[(i)]
\item $\bigE \cup \bigE^*\bigE \cup \bigE\bigE^*$ is rigid at $0$ and hyperrigid;
\item $\bigE \cup \bigE^*\bigE$ is rigid at $0$, and is hyperrigid if and only if $G$ is row-finite; and
\item $\bigE$ itself is not rigid at $0$, and is hyperrigid if and only if $G$ is row-finite.
\end{enumerate}
If, furthermore, $G$ has infinitely many vertices, then inside $\bigO(G)^1$
\begin{enumerate}[(i)]
\item[(i)$^1$] $\bigE \cup \bigE^*\bigE \cup \bigE\bigE^* \cup \{1\}$ is hyperrigid;
\item[(ii)$^1$] $\bigE \cup \bigE^*\bigE\cup\{1\}$ is hyperrigid if and only if $G$ is row-finite; and
\item[(iii)$^1$] $\bigE \cup\{1\}$ is not hyperrigid.
\end{enumerate}
\end{Tm}

\section*{Acknowledgments}
The author would like to thank Orr Shalit, his Ph.D. advisor, for his guidance and support and for some valuable comments and suggestions regarding this manuscript.
This paper continues a previous study with Adam Dor-On, and the author would like to thank him as well for a careful reading of this manuscript and some constructive comments.
A special thanks goes
to the anonymous referee, whose remarks led to an improved version
of this paper.

\end{document}